\documentclass[a4paper,11pt,AutoFakeBold]{amsart}
\usepackage{graphicx}
\usepackage{mathptmx}
\usepackage{amsmath}
\usepackage{amssymb}
\usepackage{mathrsfs}
\usepackage{enumitem}
\usepackage{xcolor}
\usepackage{cite}

\newmuskip\pFqmuskip

\newcommand*\pFq[6][8]{%
  \begingroup 
  \pFqmuskip=#1mu\relax
  \mathcode`=\string"8000
  \begingroup\lccode`\~=`\,
  \lowercase{\endgroup\let~}\pFqcomma
  F^{#2}_{#3}{\left(\genfrac..{0pt}{}{#4}{#5}\bigg|#6\right)}%
  \endgroup
}
\newcommand{\pFqcomma}{\mskip\pFqmuskip}

\newtheorem{theorem}{Theorem}[section]

\begin{document}

\title[Probabilistic Central Bell Polynomials]{Probabilistic Central Bell Polynomials}

\author{Rongrong Xu}
\address{School of Science, Xi'an Technological University, Xi'an, 710021, Shaanxi, P. R. China}
\email{xurongrong0716@163.com}
\author{Yuankui Ma}
\address{School of Science, Xi'an Technological University, Xi'an, 710021, Shaanxi, P. R. China}
\email{mayuankui@xatu.edu.cn}
\author{Taekyun  Kim*}
\address{School of Science, Xi’an Technological University, Xi’an 710021, Shaanxi, China; Department of Mathematics, Kwangwoon University, Seoul 139-701, Republic of Korea}
\email{tkkim@kw.ac.kr}
\author{Dae San  Kim }
\address{Department of Mathematics, Sogang University, Seoul 121-742, Republic of Korea}
\email{dskim@sogang.ac.kr}

\author{Salah Boulaaras }
\address{Department of Mathematics, Qassim University, Buraydah Al-Qassim, 51452, Saudi Arabia}
\email{$saleh_{-}boulaares$@yahoo.fr}

\thanks{$^*$ Corresponding authors}
\subjclass[2010]{11B73; 11B83}
\keywords{probabilistic central factorial numbers of the second kind associated with $Y$; probabilistic central Bell polynomials associated with $Y$; probabilistic central Fubini polynomials associated with $Y$}

\maketitle

\begin{abstract}
Let $Y$ be a random variable whose moment generating function exists in a neighborhood of the origin. In this paper, we study the probabilistic central Bell polynomials associated with random variable $Y$, as probabilistic extension of the central Bell polynomials. In addition, we investigate the probabilistic central factorial numbers of the second kind associated with $Y$ and the probabilistic central Fubini polynomials associated with $Y$. The aim of this paper is to derive some properties, explicit expressions, certain identities and recurrence relations for those polynomials and numbers.
\end{abstract}

\section{Introduction}
Assume that $Y$ is a random variable satisfying the moment condition (see \eqref{15-1}).
The aim of this paper is to study, as probabilistic extensions of central Bell polynomials, the probabilistic central Bell polynomials associated with $Y$, along with the probabilistic central factorial numbers of the second kind associated with $Y$ and the probabilistic central Fubini polynomials associated with $Y$. We derive some properties, explicit expressions, certain identities and recurrence relations for those polynomials and numbers. In addition, we consider the case that $Y$ is the Poisson random variable with parameters $\alpha >0$. \par
The outline of this paper is as follows. In Section 1, we recall the central factorials, the central factorial numbers of the second and the central Bell polynomials. We remind the reader of the Stirling numbers of the second kind, the Bell polynomials, the partial Bell polynomials and the complete Bell polynomials. Assume that $Y$ is a random variable such that the moment generating function of $Y$,\,\, $E[e^{tY}]=\sum_{n=0}^{\infty}\frac{t^{n}}{n!}E[Y^{n}], \quad (|t| <r)$, exists for some $r >0$. Let $(Y_{j})_{j\ge 1}$ be a sequence of mutually independent copies of the random variable $Y$, and let $S_{k}=Y_{1}+Y_{2}+\cdots+Y_{k},\,\, (k \ge 1)$,\,\, with \, $S_{0}=0$. Then we recall the probabilistic Stirling numbers of the second kind associated with $Y$, ${n \brace k}_{Y}$ and the probabilistic Bell polynomials. Section 2 is the main results of this paper.  Let $(Y_{j})_{j \ge1},\,\, S_{k},\,\, (k=0,1,\dots)$ be as in the above. We define the probabilistic central factorial numbers of the second kind associated with $Y$, $T^{Y}(n,k)$. We derive an explicit expression for $T^{Y}(n,k)$ in Theorem 2.1. We define the probabilistic central Bell polynomials associated $Y$, $B_{n}^{(c,Y)}(x)$. The generating function of $B_{n}^{(c,Y)}(x)$ is found in Theorem 2.2. Explicit expressions for $B_{n}^{(c,Y)}(x)$ are deduced in Theorems 2.3 and 2.4. An additional explicit expression for $B_{n}^{(c,Y)}(x)$ is derived in Theorem 2.15 when $Y$ is the Poisson random variable with parameter $\alpha$. We define the probabilistic central Fubini polynomials associated with $Y$, $F_{n}^{(c,Y)}(x)$. We obtain explicit expressions for $F_{n}^{(c,Y)}(x)$ in Theorems 2.5 and 2.6. As to $B_{n}^{(c,Y)}(x)$, several facts are found, namely expressions in terms of partial Bell polynomials in Theorems 2.7 and 2.10, a recurrence relation in Theorem 2.8, convolution formula in Theorem 2.9 and higher-order derivatives in Theorem 2.13. We derive some identities involving $T^{Y}(n,k)$ and the partial Bell polynomials in Theorem 2.11 and 2.12. In Theorem 2.14, we find an expression of $T^{Y}(2n,2k)$ in terms of the partial Bell polynomial. In the rest of this section, we recall the facts that are needed throughout this paper.

\vspace{0.1in}

For $n\in\mathbb{N}\cup\{0\}$, the central factorials $x^{[n]}$ are defined by
\begin{equation}
\begin{aligned}
x^{[0]}=1,\quad x^{[n]}&=x\bigg(x+\frac{n}{2}-1\bigg) \bigg(x+\frac{n}{2}-2\bigg)\cdots \bigg(x+\frac{n}{2}-n+1\bigg)\\
&=x\binom{x+\frac{n}{2}-1}{n-1}(n-1)!,\quad (n\ge 1),\quad (\mathrm{see}\ [7,10-12,17,18]).\label{1}
\end{aligned}
\end{equation}
From \eqref{1}, we note that
\begin{equation}
\bigg(\frac{t}{2}+\sqrt{\frac{t^{2}}{4}+1}\bigg)^{2x}=\sum_{n=0}^{\infty}x^{[n]}\frac{t^{n}}{n!},\quad (\mathrm{see}\ [10-12,17,18]).\label{2}
\end{equation} \par
The central factorial numbers of the second kind are defined by
\begin{equation}
x^{n}=\sum_{k=0}^{n}T(n,k)x^{[k]},\quad (n\ge 0),\quad (\mathrm{see}\ [11]).\label{3}
\end{equation}
From \eqref{2} and \eqref{3}, we have
\begin{equation}
\frac{1}{k!}\Big(e^{\frac{t}{2}}-e^{-\frac{t}{2}}\Big)^{k}=\sum_{n=k}^{\infty}T(n,k)\frac{t^{n}}{n!},\quad (k\ge 0),\quad (\mathrm{see}\ [7,11]).\label{4}
\end{equation}
The central Bell polynomials are defined by
\begin{equation}
B_{n}^{(c)}(x)=\sum_{k=0}^{n}T(n,k)x^{k},\quad (n\ge 0),\quad (\mathrm{see}\ [12]).\label{5}
\end{equation}
When $x=1$, $B_{n}^{(c)}=B_{n}^{(c)}(1)$ are called the central Bell numbers.\\
From \eqref{5}, we note that
\begin{equation}
e^{x(e^{\frac{t}{2}}-e^{-\frac{t}{2}})}=\sum_{n=0}^{\infty}B_{n}^{(c)}(x)\frac{t^{n}}{n!},\quad (\mathrm{see}\ [12,17,18]).\label{6}	
\end{equation}
From \eqref{6}, we have
\begin{equation}
B_{n}^{(c)}(x)=\sum_{l=0}^{\infty}\sum_{j=0}^{\infty}\binom{l+j}{l}(-1)^{j}\bigg(\frac{l}{2}-\frac{j}{2}\bigg)^{n}\frac{x^{l+j}}{(l+j)!},\quad (n\ge 0),\quad (\mathrm{see}\ [12]).\label{7}	
\end{equation} \par
For $n\ge 0$, the Stirling numbers of the second kind are defined as
\begin{equation}
x^{n}=\sum_{k=0}^{n}{n \brace k}(x)_{k},\quad (n\ge 0),\quad (\mathrm{see}\ [1-24]),\label{8}
\end{equation}
where $(x)_{0}=1,\ (x)_{n}=x(x-1)\cdots(x-n+1),\ (n\ge 1)$. \\
The Bell polynomials are defined by
\begin{equation}
\phi_{n}(x)=\sum_{k=0}^{n}{n\brace k}x^{k},\quad (n\ge 0),\quad (\mathrm{see}\ [7,15,16,19]).\label{9}
\end{equation}
When $x=1$, $\phi_{n}=\phi_{n}(1)$ are called the Bell numbers. \par
For any integer $k\ge 0$, the partial Bell polynomials are defined by
\begin{equation}
\frac{1}{k!}\bigg(\sum_{m=1}^{\infty}x_{m}\frac{t^{m}}{m!}\bigg)^{k}=\sum_{n=k}^{\infty}B_{n,k}\big(x_{1},x_{2},\dots,x_{n-k+1}\big)\frac{t^{n}}{n!},\quad (\mathrm{see}\ [7,13,17]).\label{10}	
\end{equation}
Thus, by \eqref{10}, we get
\begin{equation}
\begin{aligned}
&B_{n,k}(x_{1},x_{2},\dots,x_{n-k+1})\\
&=\sum_{\substack{l_{1}+\cdots+l_{n-k+1}=k\\ l_{1}+2l_{2}+\cdots+(n-k+1)l_{n-k+1}=n}} \frac{n!}{l_{1}!l_{2}!\cdots l_{n-k+1}!}\bigg(\frac{x_{1}}{1!}\bigg)^{1!} \bigg(\frac{x_{2}}{2!}\bigg)^{l_{2}}\cdots \bigg(\frac{x_{n-k+1}}{(n-k+1)!}\bigg)^{l_{n-k+1}}.
\end{aligned}\label{11}
\end{equation}
Thus, by \eqref{10}, we get
\begin{equation}
B_{n,k}(1,1,\dots,1)={n\brace k},\quad (n\ge k\ge 0).\label{12}	
\end{equation}
The complete Bell polynomials are given by
\begin{equation}
\exp\bigg(\sum_{i=1}^{\infty}x_{i}\frac{t^{i}}{i!}\bigg)=\sum_{n=0}^{\infty}B_{n}(x_{1},x_{2},\dots,x_{n})\frac{t^{n}}{n!},\quad (\mathrm{see}\ [7,13,17]).\label{13}
\end{equation}
From \eqref{10} and \eqref{13}, we note that
\begin{equation}
B_{n}(x_{1},x_{2},\dots,x_{n})=\sum_{k=0}^{n}B_{n,k}(x_{1},x_{2},\dots,x_{n-k+1}),\quad (n\ge 0).\label{14}
\end{equation}
By \eqref{9} and \eqref{14}, we get
\begin{equation}
B_{n}(x,x,\dots,x)=\phi_{n}(x),\quad (n\ge 0).\label{15}
\end{equation} \par
Assume that $Y$ is a random variable such that the moment generating function of $Y$,
\begin{equation}
E\big[e^{tY}\big]=\sum_{n=0}^{\infty}E[Y^{n}]\frac{t^{n}}{n!},\quad (|t|<r),\ \ \textrm{exists for some $r>0$}. \label{15-1}
\end{equation}
Let $(Y_{j})_{j\ge 1}$ be a sequence of mutually independent copies of the random variable $Y$, and let
\begin{displaymath}
S_{k}=Y_{1}+Y_{2}+\cdots+Y_{k},\quad (k\in\mathbb{N}),\ \mathrm{with}\quad S_{0}=0.
\end{displaymath}
The probabilistic Stirling numbers of the second kind are given by
\begin{equation}
\frac{1}{k!}\Big(E[e^{tY}]-1\Big)^{k}=\sum_{n=k}^{\infty}{n\brace k}_{Y}\frac{t^{n}}{n!},\quad (k\ge 0),\quad (\mathrm{see}\ [3,13]).\label{16}
\end{equation}
In view of \eqref{9}, the probabilistic Bell polynomials are defined by
\begin{equation}
\phi_{n}^{Y}(x)=\sum_{k=0}^{n}{n\brace k}_{Y}x^{k},\quad (n\ge 0),\quad (\mathrm{see}\ [13]).\label{17}
\end{equation}
From \eqref{16}, we have
\begin{equation*}
	{n\brace k}_{Y}=\frac{1}{k!}\sum_{j=0}^{k}\binom{k}{j}(-1)^{k-j}E\big[S_{j}^{n}\big],\quad (n\ge k\ge 0).
\end{equation*}
When $Y=1$, ${n\brace k}_{Y}={n\brace k},\quad (n\ge k\ge 0)$. \\
By \eqref{17}, we get
\begin{equation}
\sum_{n=0}^{\infty}\phi_{n}^{Y}(x)\frac{t^{n}}{n!}=e^{x(E[e^{tY}]-1)}.\label{18}	
\end{equation}
When $Y=1$, $\phi^{Y}(x)=\phi_{n}(x),\quad (n\ge 0)$.\par

\section{Probabilistic Central Bell polynomials}
Let $Y$ be a random variable and let $(Y_{j})_{j\ge 1}$ be  a sequence of mutually independent copies of the random variable $Y$ with
\begin{equation*}
S_{0}=0,\quad S_{k}=Y_{1}+Y_{2}+\cdots+Y_{k},\quad (k\ge 1).
\end{equation*}
In view of \eqref{4}, we define the {\it{probabilistic central factorial numbers of the second kind associated with $Y$}} by
\begin{equation}
\frac{1}{k!}\Big(E[e^{\frac{Y}{2}t}]-E[e^{-\frac{Y}{2}t}]\Big)^{k}=\sum_{n=k}^{\infty}T^{Y}(n,k)\frac{t^{n}}{n!},\quad (k\ge 0). \label{19}	
\end{equation}
When $Y=1$, we have $T^{Y}(n,k)=T(n,k),\ (n\ge k\ge 0)$. From \eqref{19}, we have
\begin{align}
\sum_{n=k}^{\infty}T^{Y}(n,k)\frac{t^{n}}{n!}&=\frac{1}{k!}\sum_{j=0}^{k}\binom{k}{j}\Big(E\big[e^{\frac{Y}{2}t}\big]\Big)^{k-j}(-1)^{j}\Big(E\big[e^{-\frac{Y}{2}t}\big]\Big)^{j}\label{20} \\
&=\frac{1}{k!}\sum_{j=0}^{k}\binom{k}{j}(-1)^{j}E\Big[e^{-\frac{t}{2}(Y_{1}+Y_{2}+\cdots+Y_{j})}\Big] E\Big[e^{\frac{t}{2}(Y_{1}+Y_{2}+\cdots+Y_{k-j})}\Big]\nonumber \\
&=\frac{1}{k!}\sum_{j=0}^{k}\binom{k}{j}(-1)^{j}E\Big[e^{-\frac{t}{2}S_{j}}\Big]E\Big[e^{\frac{t}{2}S_{k-j}}\Big]\nonumber \\
&=\sum_{n=0}^{\infty}\frac{1}{k!}\sum_{j=0}^{k}\sum_{l=0}^{n}\binom{n}{l}\binom{k}{j}\bigg(\frac{1}{2}\bigg)^{n}(-1)^{n-l-j}E\big[S_{k-j}^{l}\big]E\big[S_{j}^{n-l}\big]\frac{t^{n}}{n!}.\nonumber
\end{align}
Thus, by comparing the coefficients on both sides of \eqref{20}, we obtain the following theorem.
\begin{theorem}
For $n\ge k\ge 0$, we have
\begin{displaymath}
2^{n}T^{Y}(n,k)= \frac{1}{k!}\sum_{j=0}^{k}\sum_{l=0}^{n}\binom{n}{l}\binom{k}{j}(-1)^{n-l-j}E\big[S_{k-j}^{l}\big]E\big[S_{j}^{n-l}\big].
\end{displaymath}
\end{theorem}
In view of \eqref{5}, we define the {\it{probabilistic central Bell polynomials associated with $Y$}} by
\begin{equation}
B_{n}^{(c,Y)}(x)=\sum_{k=0}^{n}T^{Y}(n,k)x^{k},\quad (n\ge 0). \label{22}	
\end{equation}
When $Y=1$, $B_{n}^{(c,Y)}(x)=B_{n}^{(c)}(x),\ (n\ge 0)$. In particular, for $x=1$, $B_{n}^{(c,Y)}=B_{n}^{(c,Y)}(1)$ are called the probabilistic central Bell numbers. \\
From \eqref{22}, we note that
\begin{align}
\sum_{n=0}^{\infty}B_{n}^{(c,Y)}(x)\frac{t^{n}}{n!}&=\sum_{n=0}^{\infty}\sum_{k=0}^{n}x^{k}T^{Y}(n,k)\frac{t^{n}}{n!} \label{23} \\
&=\sum_{k=0}^{\infty}\bigg(\sum_{n=k}^{\infty}T^{Y}(n,k)\frac{t^{n}}{n!}\bigg)x^{k}\nonumber \\
&=\sum_{k=0}^{\infty}x^{k}\frac{1}{k!}\Big(E\big[e^{\frac{Y}{2}t}\big]-E\big[e^{-\frac{Y}{2}t}\big]\Big)^{k}\nonumber \\
&=e^{x(E[e^{\frac{Y}{2}t}]-E[e^{-\frac{Y}{2}t}])}.\nonumber
\end{align}
Therefore, by \eqref{23}, we obtain the following theorem.
\begin{theorem}
The generating function of the probabilistic central Bell polynomials associated with $Y$ is given by
\begin{equation}
e^{x(E[e^{\frac{Y}{2}t}]-E[e^{-\frac{Y}{2}t}])}= \sum_{n=0}^{\infty}B_{n}^{(c,Y)}(x)\frac{t^{n}}{n!}.\label{24}
\end{equation}
\end{theorem}
From \eqref{24}, we have
\begin{align}
&\sum_{n=0}^{\infty}B_{n}^{(c,Y)}(x)\frac{t^{n}}{n!}=e^{x(E[e^{\frac{Y}{2}t}])} e^{-x(E[e^{-\frac{Y}{2}t}])}	\label{25}\\
&=\sum_{l=0}^{\infty}\frac{x^{l}}{l!}E\big[e^{\frac{t}{2}S_{l}}\big]\sum_{j=0}^{\infty}\frac{(-1)^{j}}{j!}x^{j}E\big[e^{-\frac{t}{2}S_{j}}\big] \nonumber \\
&=\sum_{l=0}^{\infty}\frac{x^{l}}{l!}\sum_{m=0}^{\infty}E[S_{l}^{m}]\frac{1}{m!}\bigg(\frac{t}{2}\bigg)^{m}\sum_{j=0}^{\infty}\frac{(-1)^{j}}{j!}x^{j}\sum_{p=0}^{\infty}E\big[S_{j}^{p}\big]\frac{(-1)^{p}}{p!}\bigg(\frac{t}{2}\bigg)^{p}\nonumber \\
&=\sum_{n=0}^{\infty}\sum_{l=0}^{\infty}\sum_{j=0}^{\infty}\frac{x^{l+j}}{l!j!}(-1)^{j}\bigg(\frac{1}{2}\bigg)^{n}\sum_{m=0}^{n}\binom{n}{m}(-1)^{n-m}E\big[S_{l}^{m}\big]E\big[S_{j}^{n-m}\big]\frac{t^{n}}{n!}.\nonumber
\end{align}
Thus, by comparing the coefficients on both sides of \eqref{25}, we obtain the following theorem.
\begin{theorem}
For $n\ge 0$, we have
\begin{equation*}
B_{n}^{(c,Y)}(x)
=\frac{1}{2^{n}}\sum_{l=0}^{\infty}\sum_{j=0}^{\infty}\frac{x^{l+j}}{l!j!}(-1)^{j}\sum_{m=0}^{n}\binom{n}{m}(-1)^{n-m}E\big[S_{l}^{m}\big]E\big[S_{j}^{n-m}\big].
\end{equation*}
\end{theorem}
When $Y=1$, we get the following result in \eqref{7}:
\begin{displaymath}
B_{n}^{(c)}(x)=\sum_{l=0}^{\infty}\sum_{j=0}^{\infty}\binom{l+j}{l}(-1)^{j}\frac{x^{l+j}}{(l+j)!}\bigg(\frac{l}{2}-\frac{j}{2}\bigg)^{n}.
\end{displaymath}
Now, we observe that
\begin{align}
&\sum_{n=0}^{\infty}B_{n}^{(c,Y)}(x)\frac{t^{n}}{n!}=e^{x(E[e^{\frac{Y}{2}t}]-E[e^{-\frac{Y}{2}t}])}\label{27}\\
&=\sum_{l=0}^{\infty}x^{l}\frac{1}{l!}\Big(E\big[e^{\frac{Y}{2}t}\big]-1\Big)^{l}\sum_{j=0}^{\infty}(-1)^{j}x^{j}\frac{1}{j!}\Big(E\big[e^{-\frac{Y}{2}t}\big]-1\Big)^{j}\nonumber \\
&=\sum_{l=0}^{\infty}x^{l}\sum_{i=l}^{\infty}{i \brace l}_{Y}\bigg(\frac{1}{2}\bigg)^{i}\frac{t^{i}}{i!}\sum_{j=0}^{\infty}(-1)^{j}x^{j}\sum_{k=j}^{\infty}{k\brace j}_{Y}(-1)^{k}\bigg(\frac{1}{2}\bigg)^{k}\frac{t^{k}}{k!}\nonumber\\
&=\sum_{i=0}^{\infty}\sum_{l=0}^{i}x^{l}{i \brace l}_{Y}\bigg(\frac{1}{2}\bigg)^{i}\frac{t^{i}}{i!}\sum_{k=0}^{\infty}\sum_{j=0}^{k}(-1)^{j}x^{j}{k \brace j}_{Y}(-1)^{k}\bigg(\frac{1}{2}\bigg)^{k}\frac{t^{k}}{k!}\nonumber \\
&=\sum_{n=0}^{\infty}\bigg(\frac{1}{2}\bigg)^{n}\sum_{k=0}^{n}\sum_{j=0}^{k}\sum_{l=0}^{n-k}x^{j+l}{n-k \brace l}_{Y}\binom{n}{k}{k \brace j}_{Y}(-1)^{j+k}\frac{t^{n}}{n!}.\nonumber
\end{align}
Therefore, by \eqref{27}, we obtain the following theorem.
\begin{theorem}
For $n\ge 0$, we have
\begin{displaymath}
B_{n}^{(c,Y)}(x)= \bigg(\frac{1}{2}\bigg)^{n}\sum_{k=0}^{n}\sum_{j=0}^{k}\sum_{l=0}^{n-k}(-1)^{j+k}\binom{n}{k}{n-k \brace l}_{Y}{k \brace j}_{Y}x^{j+l}.
\end{displaymath}
\end{theorem}
Now, we define the {\it{probabilistic central Fubini polynomials associated with $Y$}} by
\begin{equation}
\sum_{n=0}^{\infty}F_{n}^{(c,Y)}(x)\frac{t^{n}}{n!}=\frac{1}{1-x(E[e^{\frac{Y}{2}t}]-E[e^{-\frac{Y}{2}t}])}.\label{28}	
\end{equation}
In particular, $F_{n}^{(c,Y)}=F_{n}^{(c,Y)}(1),\ (n\ge 0)$, are called the probabilistic central Fubini numbers. By \eqref{19} and \eqref{28}, we get
\begin{align}
\sum_{n=0}^{\infty}F_{n}^{(c,Y)}(x)\frac{t^{n}}{n!}&=\sum_{k=0}^{\infty}x^{k}k!\frac{1}{k!}\Big(E[e^{\frac{Y}{2}t}]-E[e^{-\frac{Y}{2}t}]\Big)^{k}\label{29} \\
&=\sum_{k=0}^{\infty}x^{k}k!\sum_{n=k}^{\infty}T_{n}^{Y}(n,k)\frac{t^{n}}{n!}\nonumber\\
&=\sum_{n=0}^{\infty}\sum_{k=0}^{n}x^{k}k!T_{n}^{Y}(n,k)\frac{t^{n}}{n!}.\nonumber
\end{align}
Therefore, by \eqref{29}, we obtain the following theorem.
\begin{theorem}
For $n\ge 0$, we have
\begin{displaymath}
F_{n}^{(c,Y)}(x)= \sum_{k=0}^{n}k!T_{n}^{Y}(n,k)x^{k}.
\end{displaymath}
\end{theorem}
From \eqref{28}, we note that
\begin{align}
&\sum_{n=0}^{\infty}F_{n}^{(c,Y)}(x)\frac{t^{n}}{n!}=\frac{1}{1-x(E[e^{\frac{Y}{2}t}]-E[e^{-\frac{Y}{2}t}])}\label{30} \\
&=\sum_{k=0}^{\infty}x^{k}\sum_{j=0}^{k}\binom{k}{j}\Big(E\big[e^{\frac{Y}{2}t}\big]-1\Big)^{k-j}(-1)^{j} \Big(E\big[e^{-\frac{Y}{2}t}\big]-1\Big)^{j}\nonumber \\
&=\sum_{k=0}^{\infty}x^{k}{k!}\sum_{j=0}^{k}\frac{1}{(k-j)!} \Big(E\big[e^{\frac{Y}{2}t}\big]-1\Big)^{k-j}(-1)^{j}\frac{1}{j!} \Big(E\big[e^{-\frac{Y}{2}t}\big]-1\Big)^{j}\nonumber \\
&= \sum_{k=0}^{\infty}x^{k}{k!}\sum_{j=0}^{k}(-1)^{j}\sum_{n=k}^{\infty}\sum_{i=k-j}^{n-j}{i \brace k-j}_{Y}\bigg(\frac{1}{2}\bigg)^{n}{n-i \brace j}_{Y}(-1)^{n-i}\binom{n}{i}\frac{t^{n}}{n!}\nonumber \\
&=\sum_{n=0}^{\infty}\frac{1}{2^{n}}\sum_{k=0}^{n}\sum_{j=0}^{k}\sum_{i=k-j}^{n-j}x^{k}k!{i \brace k-j}_{Y}{n-i \brace j}_{Y}\binom{n}{i}(-1)^{n-i-j}\frac{t^{n}}{n!}.\nonumber	
\end{align}
Thus, by comparing the coefficients on both sides of \eqref{30}, we obtain the following theorem.
\begin{theorem}
For $n\ge 0$, we have
\begin{displaymath}
F_{n}^{(c,Y)}(x)= \frac{1}{2^{n}}\sum_{k=0}^{n}\sum_{j=0}^{k}\sum_{i=k-j}^{n-j}(-1)^{n-i-j}k!{i \brace k-j}_{Y}{n-i \brace j}_{Y}\binom{n}{i}x^{k}.
\end{displaymath}
\end{theorem}
By \eqref{24}, we get
\begin{align}
\sum_{n=0}^{\infty}B_{n}^{(c,Y)}\frac{t^{n}}{n!}&=e^{x(                                              E[e^{\frac{Y}{2}t}]-E[e^{-\frac{Y}{2}t}])}\label{31}\\
&=\sum_{k=0}^{\infty}\frac{1}{k!}\bigg(x\sum_{j=1}^{\infty}\frac{1-(-1)^{j}}{2^{j}}E[Y^{j}]\frac{t^{j}}{j!}\bigg)^{k}\nonumber \\
&=\sum_{k=0}^{\infty}\sum_{n=k}^{\infty}B_{n,k}\Big(xE[Y],0,\frac{x}{2^{2}}E[Y^{3}],0,\dots,\frac{x(1-(-1)^{n-k+1})}{2^{n-k+1}}E[Y^{n-k+1}]\Big)\frac{t^{n}}{n!}\nonumber \\
&=\sum_{n=0}^{\infty}\sum_{k=0}^{n}B_{n,k}\Big(xE[Y],0,\frac{x}{2^{2}}E[Y^{3}],0,\dots,\frac{x(1-(-1)^{n-k+1})}{2^{n-k+1}}E[Y^{n-k+1}]\Big)\frac{t^{n}}{n!}.\nonumber
\end{align}
Therefore, by \eqref{31}, we obtain the following theorem.
\begin{theorem}
For $n\ge 0$, we have
\begin{displaymath}
B_{n}^{(c,Y)}(x)=\sum_{k=0}^{n}B_{n,k}\Big(xE[Y],0,\frac{x}{2^{2}}E[Y^{3}],0,\dots,\frac{x(1-(-1)^{n-k+1})}{2^{n-k+1}}E[Y^{n-k+1}]\Big).
\end{displaymath}
\end{theorem}
From \eqref{24}, we note that
\begin{align}
&\sum_{n=0}^{\infty}B_{n+1}^{(c,Y)}(x)\frac{t^{n}}{n!}=\frac{d}{dt}\sum_{n=0}^{\infty}B_{n}^{(c,Y)}(x)\frac{t^{n}}{n!}=\frac{d}{dt}e^{x(E[e^{\frac{Y}{2}t}]-E[e^{-\frac{Y}{2}t}])}\label{32} \\
&=x\bigg(E\bigg[\frac{Y}{2}e^{\frac{Y}{2}t}\bigg]+ E\bigg[\frac{Y}{2}e^{-\frac{Y}{2}t}\bigg]\bigg) e^{x(E[e^{\frac{Y}{2}t}]-E[e^{-\frac{Y}{2}t}])}\nonumber\\
&=x\sum_{l=0}^{\infty}\bigg(\frac{1}{2}\bigg)^{2l}E\big[Y^{2l+1}\big]\frac{t^{2l}}{(2l)!}\sum_{m=0}^{\infty}B_{m}^{(c,Y)}(x)\frac{t^{m}}{m!}\nonumber\\
&=\sum_{n=0}^{\infty}x	\sum_{l=0}^{[\frac{n}{2}]}\binom{n}{2l}\bigg(\frac{1}{2}\bigg)^{2l}E\big[Y^{2l+1}\big]B_{n-2l}^{(c,Y)}(x)\frac{t^{n}}{n!}.\nonumber
\end{align}
Therefore, by comparing the coefficients on both sides of \eqref{30}, we obtain the following theorem.
\begin{theorem}
For $n\ge 0$, we have
\begin{displaymath}
B_{n+1}^{(c,Y)}(x)= x	\sum_{l=0}^{[\frac{n}{2}]}\binom{n}{2l}\bigg(\frac{1}{2}\bigg)^{2l}E\big[Y^{2l+1}\big]B_{n-2l}^{(c,Y)}(x).
\end{displaymath}
\end{theorem}
We observe that
\begin{align}
	\sum_{n=0}^{\infty}B_{n}^{(c,Y)}(x+y)\frac{t^{n}}{n!}&=e^{(x+y)(E[e^{\frac{Y}{2}t}]-E[e^{-\frac{Y}{2}t}])}\label{33} \\
	&=e^{x(E[e^{\frac{Y}{2}t}]-E[e^{-\frac{Y}{2}t}])} e^{y(E[e^{\frac{Y}{2}t}]-E[e^{-\frac{Y}{2}t}])}\nonumber \\
	&=\sum_{k=0}^{\infty}B_{k}^{(c,Y)}(x)\frac{t^{k}}{k!}\sum_{m=0}^{\infty}B_{m}^{(c,Y)}(y)\frac{t^{m}}{m!}\nonumber \\
	&=\sum_{n=0}^{\infty}\sum_{k=0}^{n}\binom{n}{k}B_{k}^{(c,Y)}(x)B_{n-k}^{(c,Y)}(y)\frac{t^{n}}{n!}.\nonumber
\end{align}
Therefore, we obtain the following convolution formula.
\begin{theorem}[Convolution formula]
For $n\ge 0$, we have
\begin{displaymath}
B_{n}^{(c,Y)}(x+y)= \sum_{k=0}^{n}\binom{n}{k}B_{k}^{(c,Y)}(x)B_{n-k}^{(c,Y)}(y).
\end{displaymath}
\end{theorem}
By using \eqref{24}, we have
\begin{align}
	&\sum_{n=0}^{\infty}B_{n}^{(c,Y)}(x)\frac{t^{n}}{n!}= e^{x(E[e^{\frac{Y}{2}t}]-E[e^{-\frac{Y}{2}t}])}\label{34}\\
	&=\Big(e^{(E[e^{\frac{Y}{2}t}]-E[e^{-\frac{Y}{2}t}])}-1+1\Big)^{x}=\sum_{k=0}^{\infty}\binom{x}{k}k!\frac{1}{k!}\bigg(\sum_{j=1}^{\infty}B_{j}^{(c,Y)}\frac{t^{j}}{j!}\bigg)^{k}\nonumber \\
	&=\sum_{k=0}^{\infty}\binom{x}{k}k!\sum_{n=k}^{\infty}B_{n,k}\Big(B_{1}^{(c,Y)},B_{2}^{(c,Y)},\dots,B_{n-k+1}^{(c,Y)}\Big)\frac{t^{n}}{n!}\nonumber \\
	&=\sum_{n=0}^{\infty}\sum_{k=0}^{n}\binom{x}{k}k! B_{n,k}\Big(B_{1}^{(c,Y)},B_{2}^{(c,Y)},\dots,B_{n-k+1}^{(c,Y)}\Big)\frac{t^{n}}{n!}.\nonumber
\end{align}
Therefore, by \eqref{34}, we obtain the following theorem.
\begin{theorem}
For $n\ge 0$, we have
\begin{displaymath}
B_{n}^{(c,Y)}(x)= \sum_{k=0}^{n}\binom{x}{k}k! B_{n,k}\Big(B_{1}^{(c,Y)},B_{2}^{(c,Y)},\dots,B_{n-k+1}^{(c,Y)}\Big).
\end{displaymath}
\end{theorem}
It is easy to show that
\begin{equation}
t e^{x(E[e^{\frac{Y}{2}t}]-E[e^{-\frac{Y}{2}t}])}=t\sum_{j=0}^{\infty}B_{j}^{(c,Y)}(x)\frac{t^{j}}{j!}=\sum_{j=1}^{\infty}B_{j-1}^{(c,Y)}(x)j\frac{t^{j}}{j!}.\label{35}	
\end{equation}
Thus, by \eqref{35}, we get
\begin{align}
&\bigg(\sum_{j=1}^{\infty}jB_{j-1}^{(c,Y)}(x)\frac{t^{j}}{j!}\bigg)^{k}=t^{k}\Big(e^{x(E[e^{\frac{Y}{2}t}]-E[e^{-\frac{Y}{2}t}])}\Big)^{k} \label{36}\\
&=t^{k}\sum_{j=0}^{\infty}k^{j}x^{j}\frac{1}{j!}\Big(E\big[e^{\frac{Y}{2}t}\big]-E\big[e^{-\frac{Y}{2}t}\big]\Big)^{j} \nonumber \\
&=t^{k}\sum_{j=0}^{\infty}k^{j}x^{j}\sum_{n=j}^{\infty}T^{Y}(n,j)\frac{t^{n}}{n!}=\sum_{n=0}^{\infty}\sum_{j=0}^{n}k^{j}x^{j}T^{Y}(n,j)\frac{t^{n+k}}{n!}\nonumber \\
&=\sum_{n=k}^{\infty}k!\sum_{j=0}^{n-k}\binom{n}{k}k^{j}x^{j}T^{Y}(n-k,j)\frac{t^{n}}{n!}.\nonumber
\end{align}
From \eqref{36}, we note that
\begin{align}
&\sum_{n=k}^{\infty}\sum_{j=0}^{n-k}\binom{n}{k}k^{j}x^{j}T^{Y}(n-k,j)\frac{t^{n}}{n!}=\frac{1}{k!}\bigg(\sum_{j=1}^{\infty}jB_{j-1}^{(c,Y)}(x)\frac{t^{j}}{j!}\bigg)^{k}\label{37} \\
&=\sum_{n=k}^{\infty}B_{n,k}\Big(B_{0}^{(c,Y)}(x),2B_{1}^{(c,Y)}(x),3B_{2}^{(c,Y)}(x),\dots,(n-k+1)B_{n-k}^{(c,Y)}(x)\Big)\frac{t^{n}}{n!}.\nonumber
\end{align}
Therefore, by comparing the coefficients on both sides of \eqref{37}, we obtain the following theorem.
\begin{theorem}
For $n\ge k\ge 0$, we have
\begin{displaymath}
\sum_{j=0}^{n-k}\binom{n}{k}k^{j}x^{j}T^{Y}(n-k,j)= B_{n,k}\Big(B_{0}^{(c,Y)}(x),2B_{1}^{(c,Y)}(x),3B_{2}^{(c,Y)}(x),\dots,(n-k+1)B_{n-k}^{(c,Y)}(x)\Big).
\end{displaymath}
\end{theorem}
By making use of \eqref{24}, we note that
\begin{align}
&\sum_{n=k}^{\infty}B_{n,k}\Big(B_{1}^{(c,Y)}(x),B_{2}^{(c,Y)}(x),\dots,B_{n-k+1}^{(c,Y)}(x)\Big)\frac{t^{n}}{n!}=\frac{1}{k!}\bigg(\sum_{j=1}^{\infty}B_{j}^{(c,Y)}(x)\frac{t^{j}}{j!}\bigg)^{k}\label{38}\\
&=\frac{1}{k!}\Big( e^{x(E[e^{\frac{Y}{2}t}]-E[e^{-\frac{Y}{2}t}])}-1\Big)^{k}=\sum_{j=k}^{\infty}{j\brace k}x^{j}\frac{1}{j!}\Big(E[e^{\frac{Y}{2}t}]-E[e^{-\frac{Y}{2}t}]\Big)^{j} \nonumber \\
&=\sum_{j=k}^{\infty}{j \brace k}x^{j}\sum_{n=j}^{\infty}T^{Y}(n,j)\frac{t^{n}}{n!}=\sum_{n=k}^{\infty}\bigg(\sum_{j=k}^{n}{j\brace k}T^{Y}(n,j)x^{j}\bigg)\frac{t^{n}}{n!}.\nonumber
\end{align}
Therefore, by comparing the coefficients on both sides of \eqref{38}, we obtain the following theorem.
\begin{theorem}
For $n\ge k\ge 0$, we have
\begin{displaymath}
B_{n,k}\Big(B_{1}^{(c,Y)}(x),B_{2}^{(c,Y)}(x),\dots,B_{n-k+1}^{(c,Y)}(x)\Big)=\sum_{j=k}^{n}{j\brace k}T^{Y}(n,j)x^{j}.
\end{displaymath}
\end{theorem}
From \eqref{24}, we note that
\begin{align}
&\sum_{n=0}^{\infty}\bigg(\frac{d}{dx}\bigg)^{k}B_{n}^{(c,Y)}(x)\frac{t^{n}}{n!}=\bigg(\frac{d}{dx}\bigg)^{k}e^{x(E[e^{\frac{Y}{2}t}]-E[e^{-\frac{Y}{2}t}])}\label{39}\\
&=\frac{1}{k!}k!\Big(E[e^{\frac{Y}{2}t}]-E[e^{-\frac{Y}{2}t}]\Big)^{k} e^{x(E[e^{\frac{Y}{2}t}]-E[e^{-\frac{Y}{2}t}])}\nonumber \\
&=k!\sum_{l=k}^{\infty}T^{Y}(l,k)\frac{t^{l}}{l!}\sum_{j=0}^{\infty}B_{j}^{(c,Y)}(x)\frac{t^{j}}{j!}=\sum_{n=k}^{\infty}k!\sum_{j=0}^{n-k}\binom{n}{j}B_{j}^{(c,Y)}(x)T^{Y}(n-j,k)\frac{t^{n}}{n!}.\nonumber
\end{align}
Therefore, by \eqref{39}, we obtain the following theorem. \begin{theorem}
For $k\ge 1$, we have
\begin{displaymath}
\bigg(\frac{d}{dx}\bigg)^{k}B_{n}^{(c,Y)}(x)= k!\sum_{j=0}^{n-k}\binom{n}{j}B_{j}^{(c,Y)}(x)T^{Y}(n-j,k).
\end{displaymath}
\end{theorem}
From \eqref{19}, we get
\begin{align}
\sum_{n=k}^{\infty}T^{Y}(n,k)\frac{t^{n}}{n!}&=\frac{1}{k!}\Big(	E\big[e^{\frac{Y}{2}t}\big]-E\big[e^{-\frac{Y}{2}t}\big]\Big)^{k} \label{40} \\
&=\frac{1}{k!}\bigg(\sum_{j=1}^{\infty}E[Y^{j}]\bigg(\frac{1}{2}\bigg)^{j}\big(1-(-1)^{j}\big)\frac{t^{j}}{j!}\bigg)^{k}\nonumber \\
&=\sum_{n=k}^{\infty}B_{n,k}\bigg(E[Y],0,\bigg(\frac{1}{2}\bigg)^{2}E[Y^{3}],0,\dots,\bigg(\frac{1}{2}\bigg)^{n-k+1}E\big[Y^{n-k+1}\big]\big(1-(-1)^{n-k+1}\big)\bigg)\frac{t^{n}}{n!}.\nonumber
\end{align}
Thus, by \eqref{40}, we get
\begin{equation}
\begin{aligned}
	& T^{Y}(n,k)\\
	&= B_{n,k}\bigg(E[Y],0,\bigg(\frac{1}{2}\bigg)^{2}E[Y^{3}],0,\dots,\bigg(\frac{1}{2}\bigg)^{n-k+1}E\big[Y^{n-k+1}\big]\big(1-(-1)^{n-k+1}\big)\bigg),
\end{aligned}\label{41}
\end{equation}
where $n\ge k\ge 0$. \\
From \eqref{41}, we obtain the following theorem.
\begin{theorem}
	For $n\ge k\ge 0$, we have
\begin{equation*}
	\begin{aligned}
	& T^{Y}(2n,2k)\\
	&= B_{2n,2k}\bigg(E[Y],0,\bigg(\frac{1}{2}\bigg)^{2}E[Y^{3}],0,\dots,\bigg(\frac{1}{2}\bigg)^{2n-2k}E\big[Y^{2n-2k+1}\big]\bigg).
\end{aligned}
\end{equation*}
\end{theorem}
Let $Y$ be the Poisson random variable with parameter $\alpha>0$. Then we have
\begin{equation}
E\big[e^{\frac{Y}{2}t}\big]=\sum_{n=0}^{\infty}\frac{\alpha^{n}}{n!}e^{-\alpha}e^{\frac{n}{2}t}=e^{\alpha(e^{\frac{t}{2}}-1)},\label{43}	
\end{equation}
and
\begin{equation*}
E\big[e^{-\frac{Y}{2}t}\big]=\sum_{n=0}^{\infty}\frac{\alpha^{n}}{n!}e^{-\alpha}e^{-\frac{n}{2}t}=e^{\alpha(e^{-\frac{t}{2}}-1)}.
\end{equation*}
Thus, by \eqref{24} and \eqref{43}, we get
\begin{align}
&\sum_{n=0}^{\infty}B_{n}^{(c,Y)}(x)\frac{t^{n}}{n!}= e^{x(E[e^{\frac{Y}{2}t}]-E[e^{-\frac{Y}{2}t}])}\label{44}\\
&=e^{x(e^{\alpha(e^{\frac{t}{2}}-1)}-e^{\alpha(e^{-\frac{t}{2}}-1)})}\nonumber \\
&=e^{x(e^{\alpha(e^{\frac{t}{2}}-1)}-1)}e^{-x (e^{\alpha(e^{-\frac{t}{2}}-1)}-1)}\nonumber\\
&=\sum_{j=0}^{\infty}\phi_{j}(x)\frac{\alpha^{j}}{j!}\big(e^{\frac{t}{2}}-1\big)^{j}\sum_{k=0}^{\infty}\phi_{k}(-x)\frac{\alpha^{k}}{k!}\big(e^{-\frac{t}{2}}-1\big)^{k}\nonumber \\
&=\sum_{j=0}^{\infty}\phi_{j}(x)\alpha^{j}\sum_{i=j}^{\infty}{i\brace j}\bigg(\frac{1}{2}\bigg)^{i}\frac{t^{i}}{i!}\sum_{k=0}^{\infty}\phi_{k}(-x)\alpha^{k}\sum_{l=k}^{\infty}{l\brace k}(-1)^{l}\bigg(\frac{1}{2}\bigg)^{l}\frac{t^{l}}{l!}\nonumber \\
&=\sum_{i=0}^{\infty}\sum_{j=0}^{i}\phi_{j}(x)\alpha^{j}{i\brace j}\bigg(\frac{1}{2}\bigg)^{i}\frac{t^{i}}{i!}\sum_{l=0}^{\infty}\sum_{k=0}^{l}\phi_{k}(-x)\alpha^{k}{l\brace k}(-1)^{l}\bigg(\frac{1}{2}\bigg)^{l}\frac{t^{l}}{l!}\nonumber \\
&=\sum_{n=0}^{\infty}\bigg(\sum_{l=0}^{n}\sum_{k=0}^{l}\sum_{j=0}^{n-l}\phi_{j}(x)\alpha^{j+k}\binom{n-l}{j}\bigg(\frac{1}{2}\bigg)^{n}\phi_{k}(-x){l\brace k}(-1)^{l}\binom{n}{l}\bigg)\frac{t^{n}}{n!}.\nonumber
\end{align}
Therefore, by \eqref{44}, we obtain the following theorem.
\begin{theorem}
Let $Y$ be the Poisson random variable with parameter $\alpha>0$. For $n\ge 0$, we have
\begin{displaymath}
B_{n}^{(c,Y)}(x)=\frac{1}{2^{n}}\sum_{l=0}^{n}\sum_{k=0}^{l}\sum_{j=0}^{n-l} (-1)^{l}\binom{n}{l}\binom{n-l}{j}{l\brace k}\alpha^{j+k} \phi_{j}(x) \phi_{k}(-x).
\end{displaymath}
\end{theorem}
Let $Y$ be the Bernoulli random variable with probability of success $p$. Then we have
\begin{equation}
E\big[Y^{n}\big]=\sum_{i=0}^{1}i^{n}p(i)=0^{n}p(0)+1^{n}p(1)=p,\quad (n\in\mathbb{N}).\label{46}
\end{equation}
By Theorem 2.14 and \eqref{46}, we get
\begin{equation*}
\begin{aligned}
	T^{Y}(2n,2k)&=B_{2n,2k}\bigg(p,0,\bigg(\frac{1}{2}\bigg)^{2}p,0, \dots,\bigg(\frac{1}{2}\bigg)^{2n-2k}p\bigg)\\
	&=p^{2k}T(2n,2k).
\end{aligned}
\end{equation*}
\section{Conclusion}
Assume that $Y$ is a random variable such that the moment generating function of $Y$ exists in a neighborhood of the origin. In this paper, we studied by using generating functions probabilistic extensions of several special polynomials and numbers, namely the probabilistic central factorial numbers of the second kind associated with $Y$, the probabilistic central Bell polynomials associated with $Y$ and the probabilistic central Fubini polynomials associated with $Y$. \par
In more detail, we obtained an explicit expression for $T^{Y}(n,k)$ in Theorem 2.1 and a representation in terms of $B_{2n,2k}$ for $T^{Y}(2n,2k)$ in Theorem 2.14. A generating function of $B_{n}^{(c,Y)}(x)$
was derived in Theorem 2.2. We deduced explicit expressions for $B_{n}^{(c,Y)}(x)$ in Theorems 2.3 and 2.4, and in Theorem 2.15 for the special case when $Y$ is the Poisson random variable with parameter $\alpha >0$. We found explicit expressions for $F_{n}^{(c,Y)}(x)$ in Theorems 2.5 and 2.6. We deduced several facts about $B_{n}^{(c,Y)}(x)$, namely representions in terms of $B_{n,k}$ in Theorems 2.7 and 2.10, a recursive formula in Theorem 2.8, a convolution formula in Theorem 2.9 and higher-order derivatives in Theorem 2.13. We obtained identities involving $T^{Y}(n,k)$ and $B_{n,k}$ in Theorems 2.11 and 2.12. \par
It is one of our future projects to continue to study probabilistic versions of many special polynomials and numbers and to find their applications to physics, science and engineering as well as to mathematics.

\section*{Funding}
The third author of this research has been conductedby the Research Grant of Kwangwoon Universityin  2024. And this research was funded by the National Natural Science Foundation of China (No. 12271320), Key Research and Development Program of Shaanxi (No. 2023-ZDLGY-02).

\end{document}